\documentclass{ws-rv975x65}
\usepackage{ws-rv-van}             
\makeindex

\renewcommand\thesection{\arabic{section}~}

\usepackage{amsfonts}
\usepackage{amsfonts}
\usepackage{amssymb}
\usepackage{color, amsmath,amssymb, amsfonts, amstext,latexsym}

\usepackage{amssymb, epsfig, amssymb, latexsym}
\usepackage{amsmath}
\usepackage{graphicx}
\usepackage{longtable}

\numberwithin{equation}{section}

\numberwithin{theorem}{section}
\numberwithin{lemma}{section}

\numberwithin{remark}{section}
\numberwithin{example}{section}

\begin{document}

\chapter*{Quantifying Model Uncertainties\\ in the Space of Probability Measures}

\author[J. Duan, T. Gao and G. He]{Jinqiao Duan$^1$, Ting Gao$^1$ and  Guowei He$^2$ \footnote{Partly supported by   the NSF grant 1025422, NSFC grant   10971225, the open funding from the Laboratory for Nonlinear Mechanics at the Chinese Academy of Sciences, and the Fundamental Research Funds for the Central Universities, HUST 2010ZD037.}}

\address{1. Institute for Pure and Applied Mathematics\\University of California, Los Angeles, CA 90095, USA\\
\&\\Department of Applied Mathematics, Illinois Institute of Technology \\
  Chicago, IL 60616, USA  \\ \emph{E-mail: duan@iit.edu \;\;   tgao5@iit.edu } \\
\bigskip
2. Laboratory for Nonlinear Mechanics\\Institute of Mechanics, Chinese Academy of Sciences\\
 Beijing 100080, China\\  \emph{E-mail:     hgw@lnm.imech.ac.cn }    }

\newcommand{\s}{\sigma}
\newcommand{\e}{\epsilon}
\renewcommand{\a}{\alpha}
\renewcommand{\b}{\beta}
\newcommand{\om}{\omega}
\newcommand{\Om}{\Omega}
\newcommand{\D}{\Delta}
\newcommand{\p}{\partial}
\newcommand{\de}{\delta}
\newcommand{\di}{\displaystyle}
\newcommand{\N}{{\mathbb N}}
\newcommand{\R}{{\mathbb R}}

\newcommand{\F}{{\mathcal{F}}}
\newcommand{\B}{{\mathcal{B}}}
\newcommand{\cF}{\mathcal{F}}

\newcommand{\eps}{\varepsilon}
   \newcommand{\EX}{{\Bbb{E}}}
   \newcommand{\PX}{{\Bbb{P}}}

\begin{abstract}
Due to lack of scientific understanding, some mechanisms may be missing in mathematical modeling of complex phenomena in science and engineering. These mathematical models thus contain some uncertainties such as uncertain parameters. One method   to estimate these parameters is based on pathwise observations, i.e., quantifying model uncertainty in the space of sample paths for system evolution. Another method is devised here to estimate uncertain parameters, or unknown system functions,  based on experimental observations of probability distributions  for system evolution. This is called the quantification of model uncertainties in the space of probability measures. A few examples are presented to demonstrate this method, analytically or numerically.

\end{abstract}

\vskip3mm \noindent {\bf Keywords}: Model uncertainty;  Stochastic differential equations (SDEs); Probability measures; Hellinger distance; Stationary probability density; Parameter estimation; Numerical simulation

 \par

\smallskip\noindent
{\bf 2000 AMS Subject Classification : 60H40 }

\section{Introduction}  \label{intro999}

 In this chapter we   discuss some issues about quantification of model uncertainties in complex dynamical systems.


Mathematical   models for scientific and engineering   systems often involve with some uncertainties.
  We may roughly classify such uncertainties into two kinds.
  The first kind of uncertainties may be called \emph{model uncertainty}.
  They are due to  physical processes that are   not  well understood or  not well-observed, and thus are
  not or not well represented in the mathematical models.

The second kind of uncertainties may be called \emph{simulation
uncertainty}. This arises in numerical simulations of multiscale
systems that display a wide range of spatial and temporal scales,
with no clear scale separation. Due to the limitations of computer
power,    not all scales
of variability can be explicitly simulated or resolved.

These uncertainties are sometimes also called \emph{unresolved
scales}, as they are not represented (i.e.,  not  \emph{resolved})  in modeling
or simulation.  Although
these unresolved scales  may be very small or very fast, their
long time impact on the resolved simulation may be delicate (i.e.,
may be negligible or may have significant effects \cite{Crauel, XuYong}, or in other
words, uncertain). Thus, to take the effects of unresolved scales
on the resolved scales into account,   representations or
parameterizations of these effects   are desirable.

  Model uncertainties have been considered in, for example,
  \cite{Gar, Horst, Moss, Palmer2, Pas, Sura}.  Research works relevant for parameterizing unresolved scales
  include, \cite{Hasselmann, Arn00, Kantz, Stuart, Majda2, Berloff, DuanBalu, Wilks,
Williams, Sagaut, Berselli} among others. Stochastically representing unresolved scales in fluid dynamics has considered
as well \cite{Leith, Mason, Schumann}.



In this chapter, we only consider model uncertainties. Specifically, we consider dynamical systems containing uncertain parameters or unknown system functions, and examine how to estimate these parameters, using observed probability distributions of the system evolution.

After briefly comment on   estimating uncertain parameters based on observed sample paths for the system evolution in \S \ref{paths}, we then, in \S \ref{measures}, propose a method of estimating uncertain parameters based on observed probability distributions (i.e., probability measures) and present a few examples to demonstrate this method, analytically or numerically.

\section{Quantifying uncertainty in the space of paths}
\label{paths}

Since random fluctuations   are common   in the
real world, mathematical models for complex systems are often subject to uncertainties,
such as fluctuating forces, uncertain parameters, or random boundary conditions \cite{Moss, Horst, Gar, VanKampen2, WaymireDuan, Wong}.
 Stochastic differential equations (SDEs) such as
 \begin{equation}   \label{sde22222}
    dX=b(X)dt + \s(X) dB_t,
\end{equation}
 are appropriate models for many of these systems \cite{Arnold, Roz, WaymireDuan}. Here $B_t$ is a Brownian motion or Wiener process, the drift $b(X)$ and diffusion $\s(X)$   contain
 uncertain parameters (or $(b(\cdot)$ and $\s(\cdot)$ are unknown), to be estimated based on observations.

\medskip

For example,   the Langevin type models   are stochastic
differential equations describing various phenomena in physics,
biology,  and  other   fields. SDEs are used to model various price
processes, exchange rates, and interest rates, among others, in
finance. Noises in these SDEs may be modeled as a generalized time
derivative of some distinguished stochastic processes, such as
Brownian motion (BM) or other processes.



\medskip

We are interested in estimating
  parameters contained in the stochastic differential
equation \eqref{sde22222}, so that we obtain computational models useful for
investigating complex dynamics under uncertainty.

Theoretical results on parameter estimations for SDEs driven by
Brownian motion are relatively well developed
 \cite{Bishwal, Kuto, Ibragimov, Alizadeh2002, Davis2000 ,Genon-Catalot1999, Pearson1994}, and
various numerical simulations for these parameter estimations
\cite{Ait2002, Ait-Sahalia2003, Pearson1994, Nicolau2004}  are also
implemented.  See \cite{YangJ-Duan} for a more recent review about estimating and computing uncertain parameters, when dynamical systems are submit to colored or non-Gaussian noises.

 These research works on estimating uncertain parameters in dynamical systems are based on observations of sample paths.
 In the next section, we devise a method to estimate uncertain parameters based on observations of probability distributions of the system evolution.

\section{Quantifying uncertainty in the space of probability measures}
\label{measures}

Consider a dynamical system with model uncertainty, modeled by a scalar SDE
\begin{equation}\label{sde3456}
    dX=b(X)dt + \s(X) dB_t, \;\; X(0)=x_0,
\end{equation}
where the drift $b(X)$ and diffusion $\s(X)$   contain  uncertain parameters, to be estimated based on observations of probability distributions (i.e., probability measures) of the system paths $X_t$.

To this end, we need to introduce the Hellinger distance \cite{Cha} between two probability measures.  It is used    to quantify the similarity between two probability distributions.
This is a metric in the space of probability measures.

For our purpose here, we define the Hellinger distance $H(f, g)$ between two probability density functions $p(x) $ and $q(x)$ as follows
\begin{equation}\label{Hellinger}
   H^2(p, q) \triangleq  \frac12 \int_{\R^1} (\sqrt{p(x)} - \sqrt{q(x)} \; \;)^2 dx.
\end{equation}
The Hellinger distance $H(p, q)$ satisfies the property: $0 \leq H(p, q) \leq 1$.

We estimate uncertain parameters by minimizing the Hellinger distance between the true probability density $p$ for the solution process $X(t)$ and its observed probability density $q$. In reality, the  probability density $p$ has to be numerically formulated or discretized.
But in order to demonstrate the method, we consider two examples for which the true probability density $p$ can be analytically formulated.
In the first example, we minimize the Hellinger distance between the true stationary probability density  for the solution process $X(t)$ and its observed stationary probability density, while in the  second example, we do this for time-dependent probability densities.

\subsection{Observation of stationary probability distributions}

Under appropriate conditions on $b$ and $\s$ (see \cite{Klebaner}, p.170), such as, $b \leq 0$ and $\s \neq 0$ as well as some smoothness requirements, there exists a stationary probability density $p(x)$ for the SDE \eqref{sde3456}, as a solution of the steady Fokker-Planck equation $   p_{xx} + (b \sin(x)p)_x=0$,
\[ p(x)= \frac{C}{\s^2(x)} e^{\int_{x^*}^x \frac{2b(y)}{\s^2(y)}dy},  \]
where the positive normalization constant $C$ is chosen so that $p\geq 0$ and $ \int_{\R} p(x)dx=1$, i.e.,
\[C \triangleq 1/ \int_{-\infty}^{\infty} \frac{e^{\int_{x^*}^x \frac{2b(y)}{\s^2(y)}dy}}{\s^2(x)} \; dx.   \]
Note that $x^*$ here may be an arbitrary point so that the integral $\int_{x^*}^x \frac{2b(y)}{\s^2(y)}dy$ exists (say, take $x^*=0$ if that is possible).

\begin{example}
(i) A special case:  Langevin equation\\
\[ dX = -bXdt + dB_t,   \]
with parameter $b>0$. Given an ``observation" of the stationary probability density $q(x)= \frac{1}{\sqrt{\pi}} e^{- x^2}$.
Find a $b$ so that the Hellinger distance $F(b) = \frac12 \int_{\R} [\sqrt{p(x)} -\sqrt{q(x)}]^2 dx$ is minimized.

(ii) A more general case: \\
\[  dX=b(X)dt +  dB_t,  \]
with function $b(x) \leq 0$. Given an ``observation" of the stationary probability density $q(x)= \frac{1}{ \pi} \frac1{1+x^2}$ (the Cauchy  distribution).
Find a   function $b(x) \leq 0$ so that the Hellinger distance $F(b(x)) = \frac12 \int_{\R} [\sqrt{p(x)} -\sqrt{q(x)}]^2 dx$ is minimized.
\end{example}

{\bf Solution:}\\
(i)  The true stationary probability density   for the solution process $X_t$   is
$$p(x)=\frac{\sqrt{b}}{\sqrt{\pi}}e^{-bx^2}.$$

 Insert $p, q$ into the Hellinger distance $F(b)$, which is now an algebraic function of parameter $b>0$.    Thus we use deterministic calculus to find a minimizer $b$ (possibly by hand, or Matlab if needed). Note: $\int_{\R} e^{-z^2} dz =\sqrt{\pi}$.

To minimize the Hellinger distance $F(b)$, we calculate its derivative
\begin{align*}
F'(b)&=\frac{\di 1}{\di 2}\int_{\mathbb{R}}e^{-bx^2}(\frac{\di
1}{\di 2\sqrt{\pi b}}-\frac{\di \sqrt{b}}{\di \sqrt{\pi}}
x^2)dx-\frac{\di 1}{\di \sqrt{\pi}}\di\int_{\mathbb{R}}
e^{-\frac{(1+b)x^2}{2}}( \frac{1}{4} b^{-\frac{3}{4}} -\frac{x^2}{2}
b^{\frac{1}{4}})dx\\& =\frac{\di 1}{\di
\sqrt{2}}b^{\frac{1}{4}}(1+b)^{-\frac{3}{2}}-\frac{1}{2\sqrt{2}}b^{-\frac{3}{4}}(1+b)^{-\frac{1}{2}}=0.
\end{align*}

Therefore,
$$b^{\frac{1}{4}}(1+b)^{-\frac{3}{2}}=\frac{b^{-\frac{3}{4}}(1+b)^{-\frac{1}{2}}}{2}.$$
Thus we obtain the parameter $b=1$.

(ii) The true stationary probability density for the solution process $X_t$ is
$$p(x)=\frac{e^{2\int_{0}^xb(y)dy}}{\int_{-\infty}^{\infty} e^{2\int_{0}^x b(y)dy dx}}.$$

 Insert $p, q$ into the Hellinger distance  $F(b(x))$, which  is now a  functional of   $b(x)$ and thus we use   calculus of variations (on $F(b(x))$) to find a minimizer $b(x)$. We then derive the Euler-Lagrange equation to be satisfied by $b(x)$, together with appropriate boundary conditions (needed for $p(x) \geq 0$ and $\int_{\R} p(x) \; dx=1$).

 To this end, we calculate, for an arbitrary ``variations" $h(x)$
\begin{align*}
&F(b(x)+\varepsilon h(x))\\
&=\frac{1}{2} \int_\mathbb{R} \left[ \frac{e^{2\int_{0}^x(b(y)+
\varepsilon h(y)) dy }}{\int_{-\infty}^{\infty}
e^{2\int_0^x(b(y)+\varepsilon h(y)) dy} dx} -
\frac{2e^{\int_0^x(b(y)+\varepsilon h(y)) dx} }{\sqrt{\pi}
\sqrt{1+x^2}
\sqrt{\int_{-\infty}^{\infty}e^{2\int_0^x(b(y)+\varepsilon h(y))
dy}dx}}
+\frac{1}{\pi(1+x^2)} \right]dx.
\end{align*}

The Euler-Lagrange equation for $b(x)$ comes from:  $\frac{d}{d\varepsilon} |_{\varepsilon=0} \; F(b(x)+\varepsilon h(x)) =0$ for arbitrary
``variations" $h(x)$.
 In fact, the Euler-Lagrange equation for $b(x)$ is
  $$ \int_{\mathbb{R}}  e^{\int_0^x b(w) dw}
\left( e^{\int_0^x b(w) dw} - \frac{\sqrt{\int_\mathbb{R}e^{2
\int_0^x b(w)dw} dx }}{\sqrt{\pi(1+x^2)}} \right) \cdot
\int_{\mathbb{R}} \left(e^{2\int_0^z b(w)dw} \int_z^x h(y)dy \right)
dz=0.$$

After changing the order of integration (first on $y$ and $z$ then on
$y$ and $x$), we have

$$ \left\{ \int_y^{\infty} e^{\int_0^x b(w) dw} \left( e^{\int_0^x b(w) dw}  - \frac{\sqrt{\int_{\mathbb{R}}
e^{2\int_0^x b(w)dw }dx}}{\sqrt{\pi(1+x^2)}} \right)dx \cdot
\int_{-\infty}^y e^{2\int_0^zb(w)dw} dz \right\} h(y) dy=0 $$ holds
for all $h(y)$.

Therefore,
$$\int_y^{\infty} e^{\int_0^x b(w) dw} \left( e^{\int_0^x b(w) dw}  - \frac{\sqrt{\int_{\mathbb{R}}
e^{2\int_0^x b(w)dw }dx}}{\sqrt{\pi(1+x^2)}} \right)dx  \cdot
\int_{-\infty}^y e^{2\int_0^zb(w)dw} dz =0.$$

Since $\int_{-\infty}^y e^{2\int_0^zb(w)dw} dz >0$, we further obtain
$$\int_y^{\infty} e^{\int_0^x b(w) dw} \left( e^{\int_0^x b(w) dw}  - \frac{\sqrt{\int_{\mathbb{R}}
e^{2\int_0^x b(w)dw }dx}}{\sqrt{\pi(1+x^2)}} \right)dx =0$$ for
$y\in R.$
Then taking the derivative of the above equation with
respect to $y$, we arrive at
 $$e^{\int_0^y b(w)dw} =
\frac{\sqrt{\int_{\mathbb{R}} e^{2\int_0^x b(w)dw
}dx}}{\sqrt{\pi(1+y^2)}}, \ \ \ \forall y \in \mathbb{R}^1.$$

Thus, after taking the `square' and $`\ln'$ on both sides of the
above equation, we get

$$2\int_0^y b(w)dw = \ln(\int_{\mathbb{R}} e^{2\int_0^x b(w)dw
}dx)-\ln(\pi(1+y^2)).$$
Finally taking the derivative with respect to $y$, we have
$$
b(y) =
-\frac{y}{\pi(1+y^2)}\ \ \ y\in \mathbb{R}^1.
$$

Also note that we only need $\int_0^y b(w) dw <0$ for all $y\in
\mathbb{R}^1$ for the stationary probability density to make sense.\\\\\\

\subsection{Observation of   time-dependent probability distributions  }

Consider a scalar SDE
\begin{equation}\label{sde34567}
    dX=b(X)dt + \s(X) dB_t, \; X(0)=x_0.
\end{equation}
The Fokker-Planck equation \cite{Oksendal, VanKampen3, Gardiner}  for the probability density $p(x, t) \triangleq p(x, t; x_0, 0)$ for the solution $X(t, x_0)$ is
\begin{equation}\label{fpe}
    p_t = \frac12 ( \s^2(x) p(x, t))_{xx} - (b(x)p(x, t))_x,  \;\; p(x, 0)=\delta(x_0).
\end{equation}
With an observation of $p(x, t)$, we can estimate parameters, or $b(\cdot)$, or $\s(\cdot)$, by examining the inverse problem of the
Fokker-Planck equation \eqref{fpe}. For more information about inverse problems of partial differential equations, see \cite{Isakov}.

Let us look at a specific example.
\begin{example}
Consider a scalar SDE
\begin{equation}\label{sde345678}
    dX=-b \sin(X)dt + \sqrt{2} \; dB_t, \; X(0)=0.
\end{equation}
(i) Assume that an observation obtained for $p$ to be
 \begin{equation}\label{pdf}
q_1(x, t) = \frac1{2 \sqrt{\pi t}} \; e^{-\frac{x^2}{4t}}.
\end{equation}
Find the parameter $b$ by minimizing the Hellinger distance $H(p, q_1)$.

(ii) Assume that another observation obtained for $p$ to be
 \begin{equation}\label{pdf2}
 q_2(x,t)=\frac{\sqrt{t}}{\pi(t+x^2)}.
\end{equation}
Find the parameter $b$ by minimizing the Hellinger distance $H(p, q_2)$.

{\bf Solution:}\\
The Fokker-Planck equation for \eqref{sde345678} is
\begin{equation}\label{fpe2}
    p_t =   p_{xx} + (b \sin(x)p)_x, \;\; p(x, 0)=\delta(0).
\end{equation}
 In this case, we define the Hellinger distance:
$$H_i^2(b)=\max_{t\in [0,T]}\int_{-\infty}^{\infty}
(\sqrt{q_i(x,t)}-\sqrt{p(x,t,b)})^2 dx$$
where $i=1,2$, and $T$ is the time period when $q_i(x, t)$ are observed. We numerical find
  $b$ by minimizing $H_i(b)$.


The observation $q_1(x,t)$ is plotted in Figure \ref{q1plot}.
 \begin{figure}[h]
 \begin{center}
\includegraphics*[width=6.5cm,height=5cm]{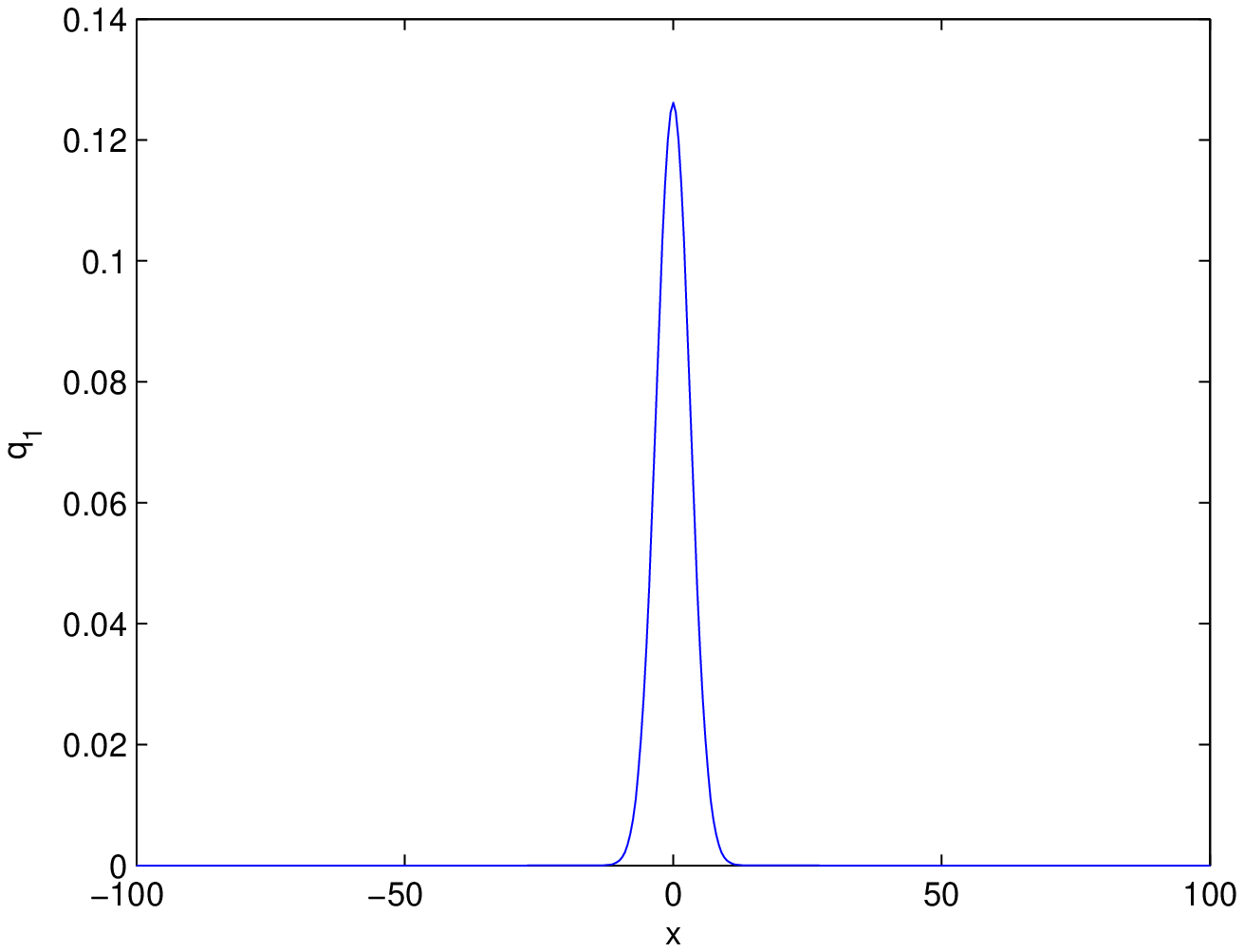}
\includegraphics*[width=6.5cm,height=5cm]{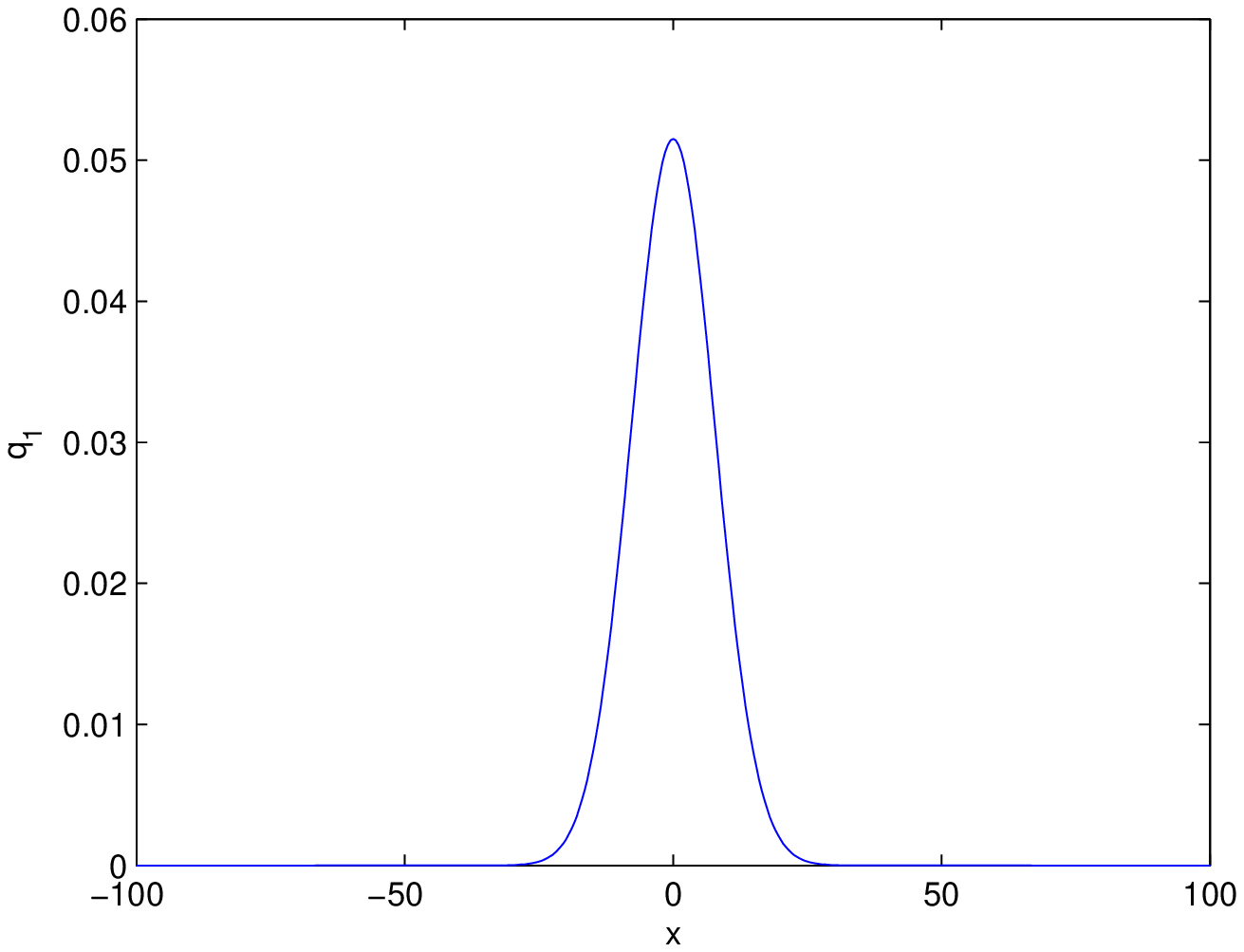}
\end{center}
 \caption{Observation   $q_1(x,t)$:  $t=5$ (top) and $t=30$ (bottom).}\label{q1plot}
\end{figure}

\newpage

 By the definition $H_1^2(b)=\max_{t\in
[0,T]}\int_{-\infty}^{\infty} (\sqrt{q_1(x,t)}-\sqrt{p(x,t,b)})^2
dx$, we have the plot of $H_1(b)$ in Figure \ref{pq1}. And
whatever $T$ is , $H_1(b)$ is always minimized when $b=0$. This gives us the parameter value $b=0$.

 \begin{figure}[h]
 \begin{center}
\includegraphics*[width=6.5cm,height=5cm]{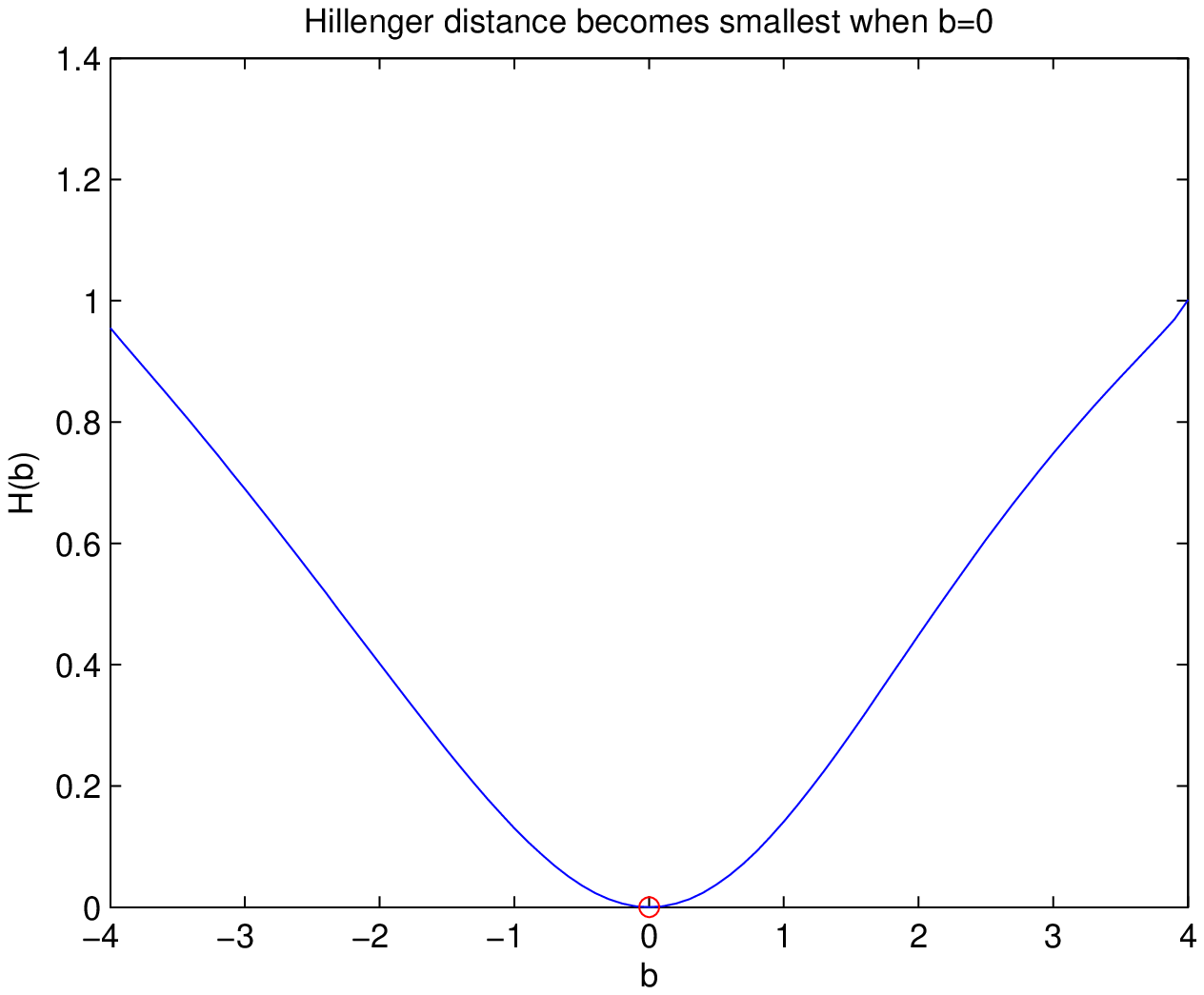}
\includegraphics*[width=6.5cm,height=5cm]{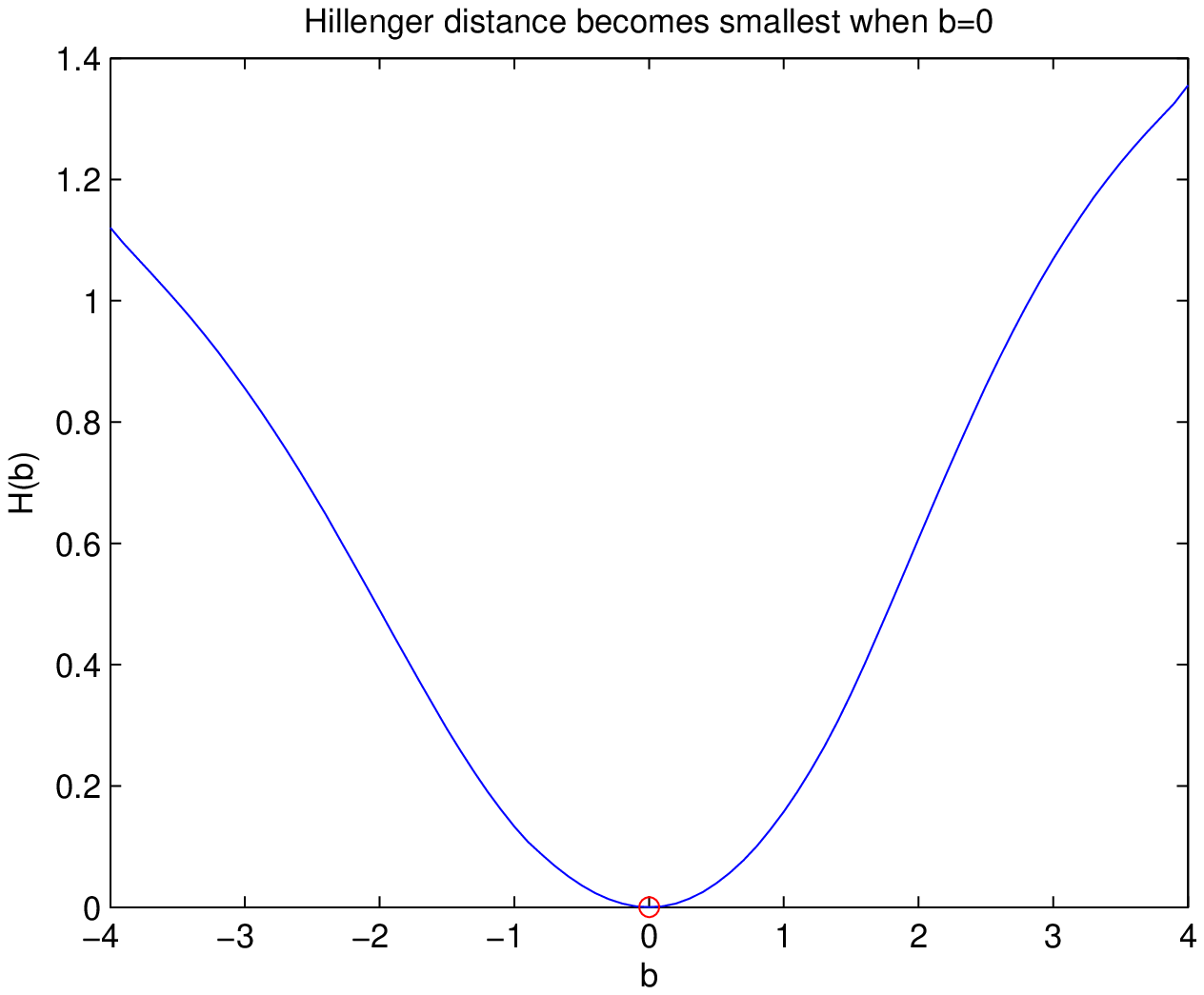}
\end{center}
 \caption{Squared Hellinger distance  between   $p(x,t)$ and the observation $q_1(x,t)$: $T=5$ (top) and $T=30$ (bottom).} \label{pq1}
\end{figure}

The observation $q_2(x,t)$ is plotted in Figure \ref{q2plot}.
 \begin{figure}[h]
 \begin{center}
\includegraphics*[width=6.5cm,height=5cm]{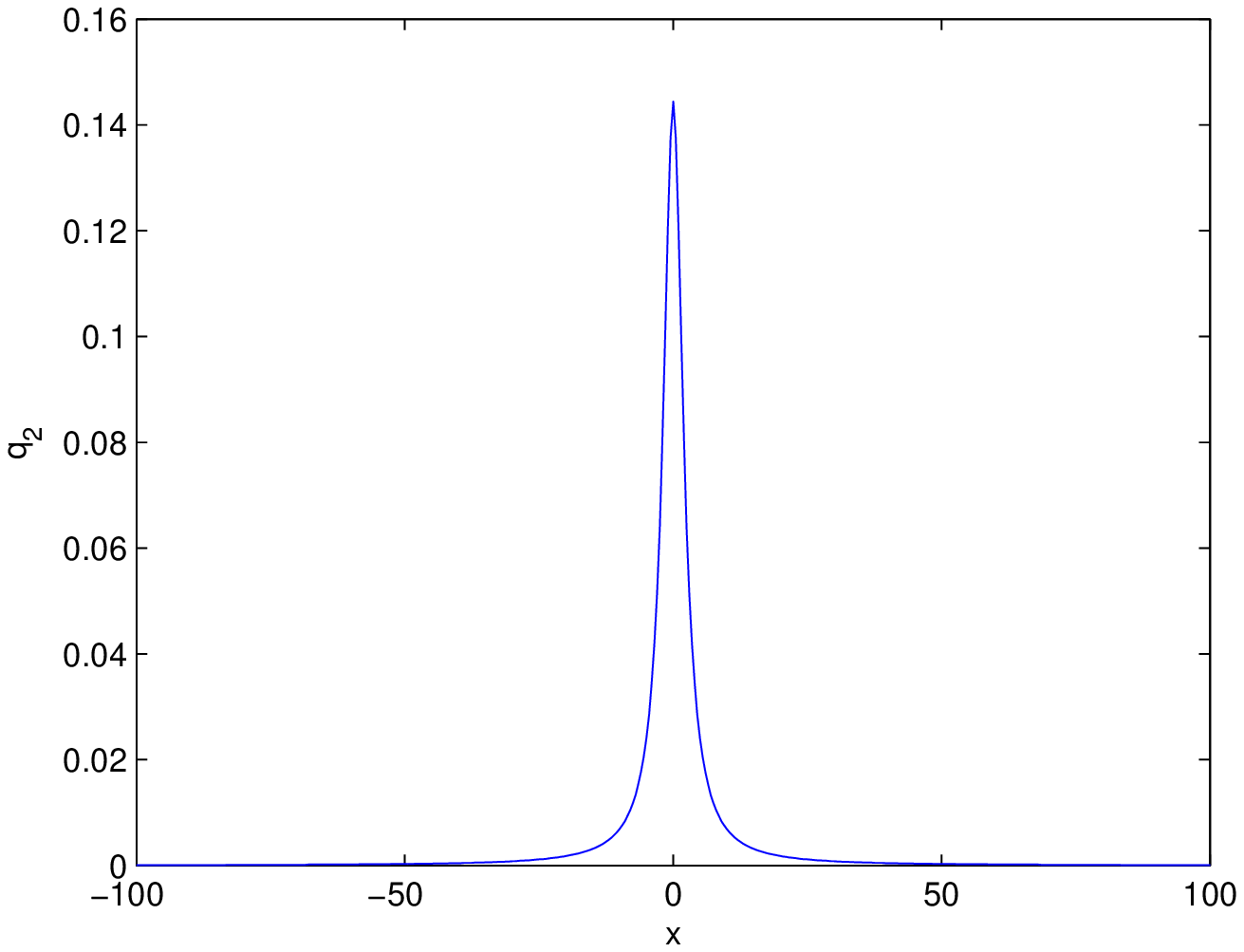}
\includegraphics*[width=6.5cm,height=5cm]{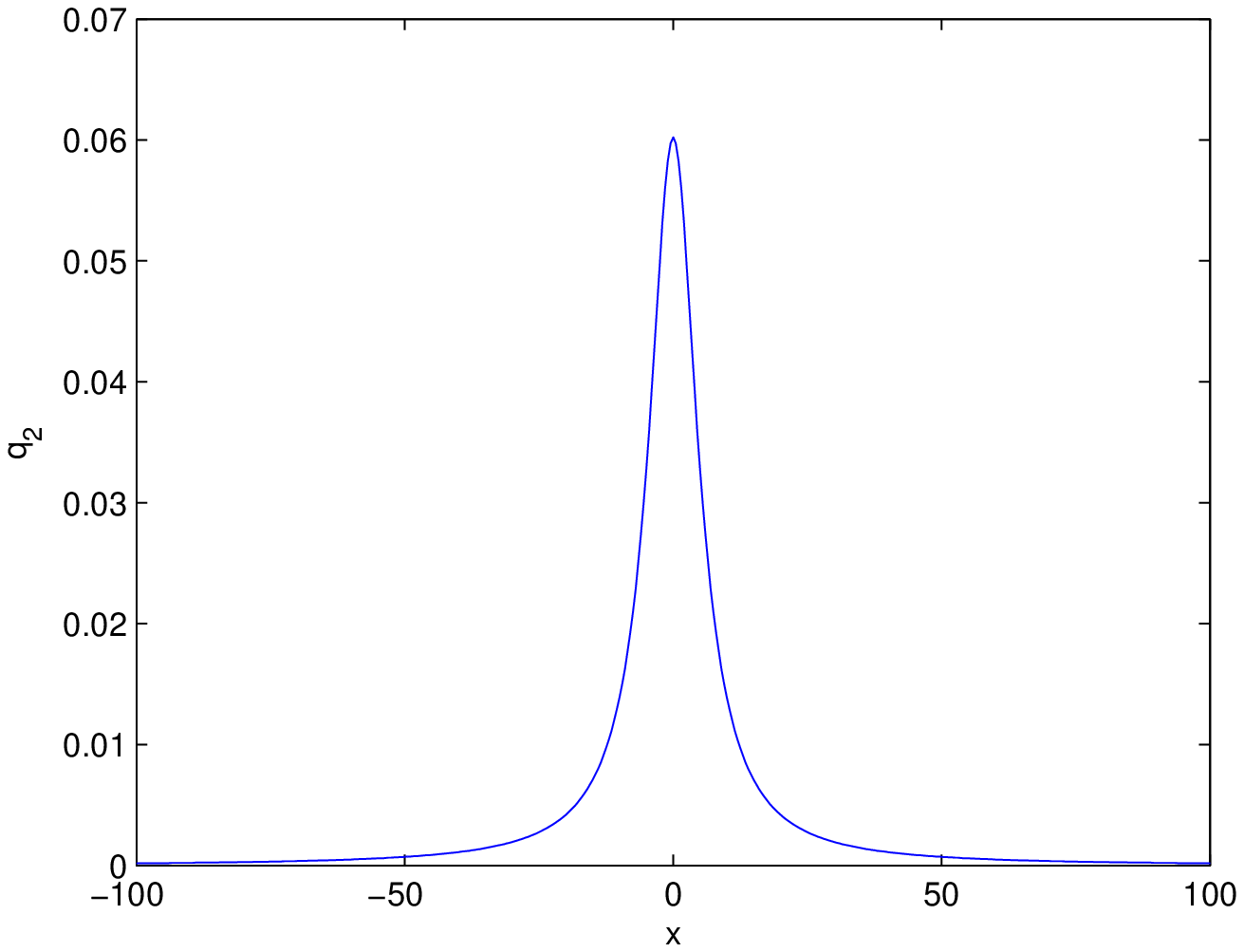}
\end{center}
 \caption{Observation   $q_2(x,t)$: $t=5$ (top) and $t=30$ (bottom).} \label{q2plot}
\end{figure}

By the definition $H_2^2(b)=\max_{t\in [0,T]}\int_{-\infty}^{\infty}
(\sqrt{q_2(x,t)}-\sqrt{p(x,t,b)})^2 dx$, we have the plot of
$H_2(b)$ in Figure \ref{pq2}. So we   see that if $T=5$,
$H_2(b)$ is minimized when $b=0.7$ and if $T=30$, $H_2(b)$ is
minimized when $b=0.6$.

 \begin{figure}[h]
 \begin{center}
\includegraphics*[width=6.5cm,height=5cm]{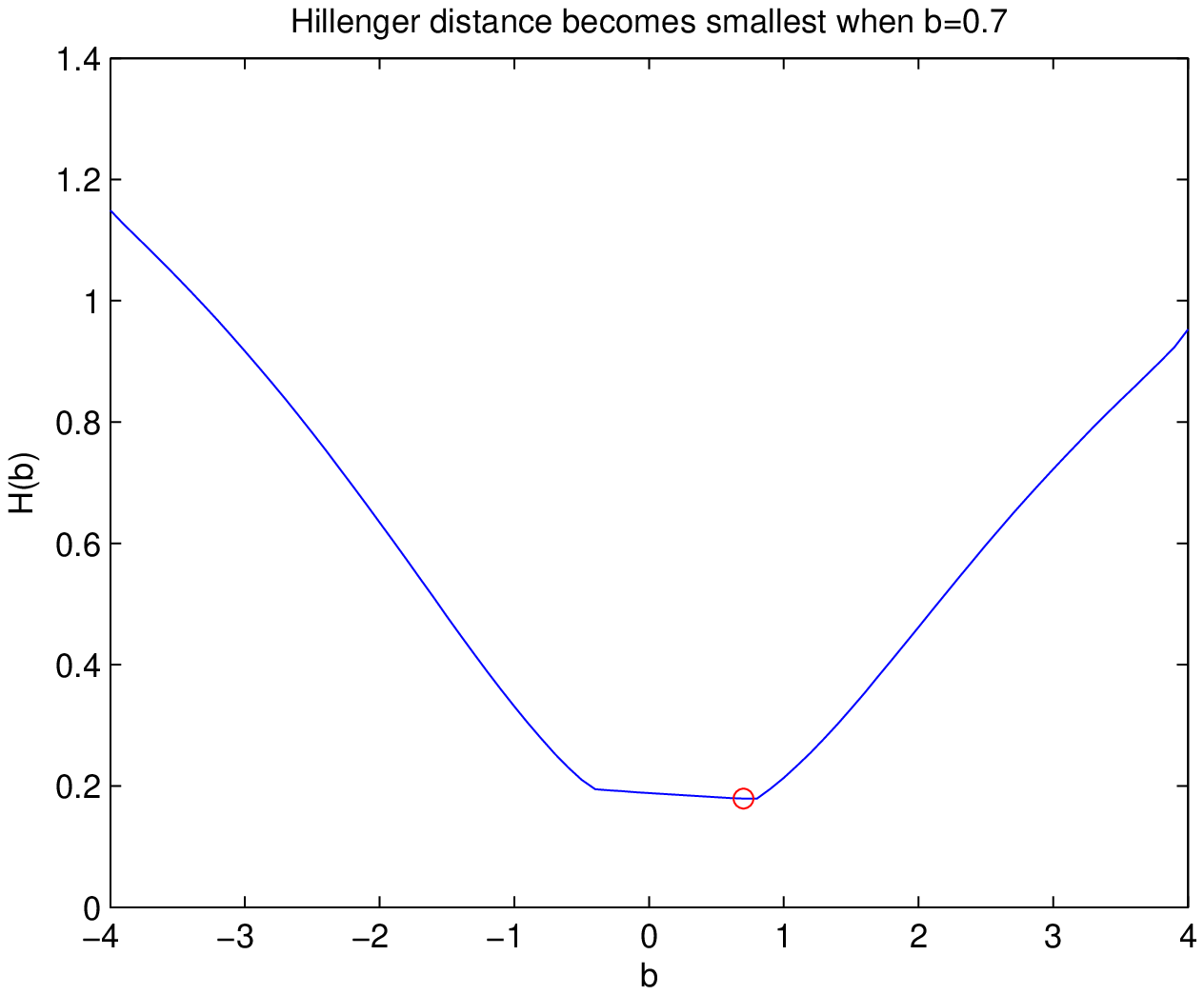}
\includegraphics*[width=6.5cm,height=5cm]{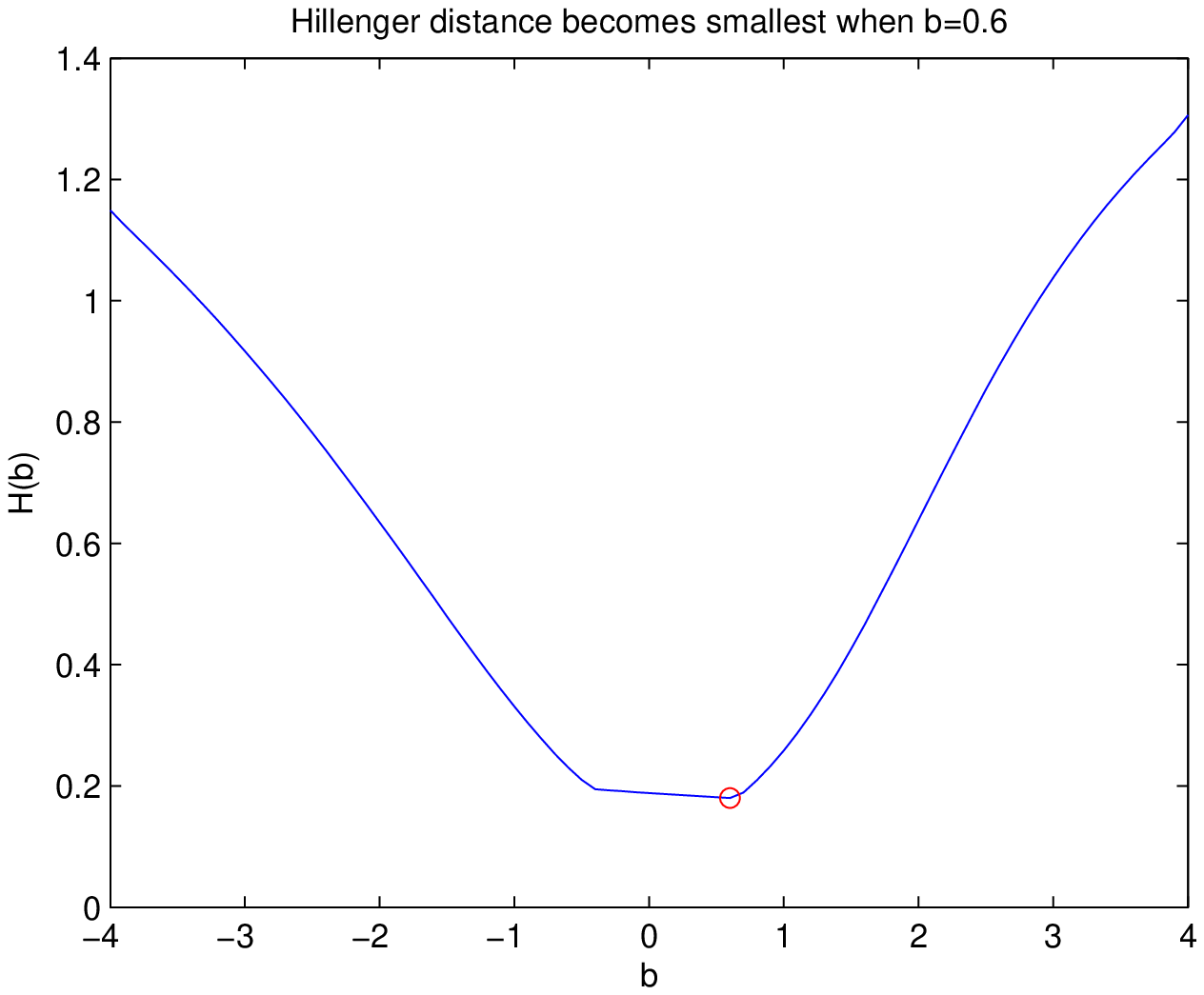}
\end{center}
 \caption{Squared Hellinger distance between   $p(x,t)$ and the observation $q_2(x,t)$: $T=5$ (top) and $T=30$ (top).} \label{pq2}
\end{figure}

\end{example}

\bigskip

\noindent {\bf Acknowledgements.} We would like to thank Huijie Qiao and Xiangjun Wang for helpful discussions.



\end{document}